\documentclass[a4paper]{article}
\usepackage{amsmath,amssymb,amsthm,amscd,euscript}
\allowdisplaybreaks  \def\lb{\label}

\def\bc{\begin{center}}
\def\ec{\end{center}}
\DeclareMathOperator{\Sine}{Sine}
\DeclareMathOperator{\SL}{SL}

 \def\pd#1#2{\frac{\partial{#1}}{\partial{#2}}}
\def\pd1#1{\frac{\partial}{\partial#1}}
\def\d1#1{\frac{d}{d#1}}

\begin{document} \author{Alexey V. Shchepetilov\\ Department of Physics, Moscow State University,
\\ 119992 Moscow, Russia\footnote{email: alexey@quant.phys.msu.su}}

\title{The geometric sense of R. Sasaki connection}
\date{}
\maketitle
\begin{abstract}
For the Riemannian manifold $M^{n}$ two special connections on the
sum of the tangent bundle $TM^{n}$ and the trivial one-dimensional
bundle are constructed. These connections are flat if and only if
the space $M^{n}$ has a constant sectional curvature $\pm 1$. The
geometric explanation of this property is given. This construction
gives a coordinate free many-dimensional generalization of the
connection from the paper: R.~Sasaki 1979 Soliton equations and
pseudospherical surfaces, Nuclear Phys., {\bf 154 B}, pp. 343-357.
It is shown that these connections are in close relation with the
imbedding of $M^{n}$ into Euclidean or pseudoeuclidean
$(n+1)$-dimension spaces.
\end{abstract}
Keywords: Riemannian space, bundle, connection, curvature,
zero-curvature representation.

PACS number: 02.40.Ky

Mathematical Subject Classification: primary 53C07; secondary
53C42, 35Q53.

\newpage
\section{Introduction} \label{Introduction}\markright{\ref{Introduction}
Introduction} In the paper \cite{Sasa} there was proposed the
formula for some local connection on a 2-dimensional real
Riemannian manifold $M^{2}$, which has played a big role in the
theory of nonlinear integrable partial differential equations. The
construction of this connection is as follows\footnote{The
following description of the Sasaki construction is slightly
different from the original one for the better agreement with the
sequel. Particulary we choose the sign in the null curvature
condition in more geometric way.}.

Let $g$ be a Riemannian metric on $M^{2}$, $\nabla$ the
corresponding Levi-Civita connection on $TM^{2}$,
$\{e_{1},e_{2}\}$ be a moving orthonormal frame on some open
domain $U\subset M^{2}$ and $\{\omega^{1},\omega^{2}\}$ a
corresponding moving coframe. The relations
$\nabla(e_{i})=\omega_{i}^{j}\otimes e_{j}$ define the connection
1-form matrix $\omega_{i}^{j}$ with respect to the frame
$\{e_{1},e_{2}\}$. The orthonormality of this frame implies that
$\omega_{1}^{1}=\omega_{2}^{2}=0,\,\omega_{1}^{2}=-\omega_{2}^{1}$.
The Levi-Civita connection is torsion-free that gives the
following structural equations:
\begin{equation} \lb{StrEq}
d\omega^{1}=\omega^{2}\wedge\omega_{2}^{1},\;d\omega^{2}=\omega^{1}\wedge\omega_{1}^{2}.
\end{equation}
The Gaussian curvature $K$ of the space $M^{2}$ is defined by
\begin{equation} \lb{GausCur}
d\omega^{1}_{2}=K\omega^{1}\wedge\omega^{2}.
\end{equation}
R. Sasaki proposed to consider a matrix 1-form on $M^{2}$:
\begin{equation} \lb{NewConnect}
A=\frac12\left(\begin{array}{cc} \omega^{2} &
-\omega^{1}+\omega^{1}_{2} \\ -\omega^{1}-\omega^{1}_{2} &
-\omega^{2} \end{array}\right),
\end{equation}
as a new connection form for some (non-specified) bundle over
$M^{2}$. The key property of the matrix 1-form $A$ is that it
satisfies the null curvature condition $\Omega\equiv dA+A\wedge
A\equiv 0$ iff $K\equiv -1$ on $U$.

In some preceding (for example \cite{Cram}) and many consequent
papers (some of the most recent are \cite{Reyes} -- \cite{Mar})
there were discussed different matrix 1-forms, depending on a
function u (or some functions) of some independent variables, such
that the null curvature condition for this form is equivalent to
one of the well known nonlinear partial differential equations
(KdV, mKdV, sine-Gordon, sinh-Gordon, non-linear Schr\"{o}dinger,
Burgers) possessing many conservation laws and reach symmetry
groups. R. Sasaki was the first who connect the matrix 1-form $A$
with the surfaces of constant negative curvature. In the paper
\cite{ChTe} there was defined a class of differential equations
$\mathcal{F}[u]=0$, which could be obtained as a null curvature
condition for the form (\ref{NewConnect}), depending on the
function $u$.

However, it seems that the geometric meaning of the connection
$\widetilde{\nabla}^{h}$, corresponding to the matrix 1-form $A$,
remained unclear. Firstly, the definition (\ref{NewConnect}) is
valid only for a local trivialization of potential unknown bundle,
because it might happen that there is no global moving frame on
$M^{2}$, for example for $M^{2}= \mathbb{S}^{2}$. Secondly,
according to (\ref{NewConnect}) the matrix 1-form $A$ is a $
\mathfrak{sl}(2, \mathbb{R})$-valued 1-form. It seems a bit
strange, because it is defined for an arbitrary metric $g$, and
the corresponding group $\SL(2,\mathbb{R})$ is the isometry group
only for $M^{2}$ equal to the hyperbolic plane $ \mathbb{H}^{2}$.
Thirdly, the forms $\omega^{1},\omega^{2}$ and $\omega^{1}_{2}$
play the similar role in (\ref{NewConnect}), but their geometric
sense is quite different. The Levi-Civita connection 1-form for
the tangent bundle $TM^{2}$ with respect to the moving frame
$\{e_{1},e_{2}\}$ is
\begin{equation}
\left(\begin{array}{cc} 0 & \omega^{1}_{2} \\ -\omega^{1}_{2} & 0
\end{array}\right).
\end{equation}
It is contained in $A$ with the strange factor $\frac12$,
violating the geometric sense of this summand due to nonlinear
relation between a connection 1-form $A$ and a curvature 2-form
$\Omega$. At last it seems to be unclear how to generalize this
construction for higher dimensions. Below we will answer on all
these questions.

Note the difference of R. Sasaki connection from Sasakian geometry
introduced in \cite{SS}.

\section{Reformulation of Sasaki's construction}
\label{Reformulation}\markright{\ref{Reformulation} Reformulation
of Sasaki's construction} Denote
$$
\sigma_{1}=\frac12\left(\begin{array}{cc} 0 & -1 \\ -1 & 0
\end{array}\right),\,
\sigma_{2}=\frac12\left(\begin{array}{cc} 1 & 0 \\ 0 & -1
\end{array}\right),\,
\sigma_{3}=\frac12\left(\begin{array}{cc} 0 & 1 \\ -1 & 0
\end{array}\right)
$$
the base in the Lie algebra $ \mathfrak{sl}(2, \mathbb{R})$. The
commutative relation for it are:
\begin{equation}\lb{commute}[\sigma_{1},\sigma_{2}]=\sigma_{3},\,[\sigma_{2},\sigma_{3}]=-\sigma_{1},
\,[\sigma_{3},\sigma_{1}]=-\sigma_{2}.\end{equation}
Then the
connection matrix 1-form (\ref{NewConnect}) can be expressed as
\begin{equation}\lb{NewForm}
A=\sigma_{1}\omega^{1}+\sigma_{2}\omega^{2}+\sigma_{3}\omega^{1}_{2}
\end{equation}
and the corresponding curvature form will be
\begin{align*}
\Omega&=A+\frac12[A,A]=\sigma_{1}d\omega^{1}+\sigma_{2}d\omega^{2}+\sigma_{3}d\omega^{1}_{2}+
[\sigma_{1},\sigma_{2}]\omega^{1}\wedge\omega^{2}\\ &+
[\sigma_{1},\sigma_{3}]\omega^{1}\wedge\omega^{1}_{2}+
[\sigma_{2},\sigma_{3}]\omega^{2}\wedge\omega^{1}_{2}.
\end{align*}
Here we used a standard notation
$[B,C]=\sum_{i,j}[B_{i},C_{j}]\omega^{i}\wedge\omega^{j}$, where
$B=\sum_{i}B_{i}\omega^{i},\,C=\sum_{i}C_{i}\omega^{i};\;B_{i},C_{i}$
are coefficients in Lie algebra $ \mathfrak{g}$ and $\omega^{i}$
are scalar differential 1-forms. When $ \mathfrak{g}$ is a matrix
algebra it is obvious that $B\wedge C=\dfrac12[B,C]$. Hence the
condition $\Omega=0$ depends only on relations (\ref{StrEq}),
(\ref{GausCur}) and commutative relations in the algebra
$\mathfrak{sl}(2,\mathbb{R})$.

It is well known that Lie algebras $ \mathfrak{sl}(2,
\mathbb{R}),\mathfrak{so}(2,1)$ and $\mathfrak{su}(1,1)$ are
isomorphic, so we can change elements $
\sigma_{1},\sigma_{2},\sigma_{3}$ in (\ref{NewForm}) by the
equivalent elements from $\mathfrak{so}(2,1)$:
$$
\bar\sigma_{1}=\left(\begin{array}{ccc} 0 & 0 & 1 \\ 0 & 0 & 0 \\
1 & 0 & 0 \end{array}\right),\,
\bar\sigma_{2}=\left(\begin{array}{ccc} 0 & 0 & 0 \\ 0 & 0 & 1 \\
0 & 1 & 0 \end{array}\right),\,
\bar\sigma_{3}=\left(\begin{array}{ccc} 0 & 1 & 0 \\ -1 & 0 & 0 \\
0 & 0 & 0 \end{array}\right)
$$
with the same commutative relations (\ref{commute}). After this
substitution the 1-form $A$ becomes:
\begin{equation}\lb{newform}
A=\bar\sigma_{1}\omega^{1}+\bar\sigma_{2}\omega^{2}+\bar\sigma_{3}\omega^{1}_{2}=
\left(\begin{array}{ccc} 0 & \omega_{2}^{1} & \omega^{1} \\
-\omega_{2}^{1} & 0 & \omega^{2} \\ \omega^{1} & \omega^{2} & 0
\end{array}\right),
\end{equation}
where expression
$$
\bar\sigma_{3}\omega^{1}_{2}=\left(\begin{array}{ccc} 0 &
\omega_{2}^{1} & 0 \\ -\omega_{2}^{1} & 0 & 0 \\ 0 & 0 & 0
\end{array}\right)
$$
contains the Levi-Civita connection form
\begin{equation}
\left(\begin{array}{cc} 0 & \omega^{1}_{2} \\ -\omega^{1}_{2} & 0
\end{array}\right)
\end{equation}
as a direct summand. Due to the last fact it is now possible to
give a geometric interpretation of the Sasaki's connection.

Let $E=\mathbb{R}\times M^{2}$ be a trivial one dimensional vector
bundle over $M^{2}$ and $F=TM^{2}\oplus E$ be a direct sum of two
bundles over the same base. Define an indefinite metric
$\widetilde{g}_{h}$ on each fiber $T_{x}M^{2}\oplus\mathbb{R}$ of
$F$ as a direct sum of the metric $g$ and the metric
$(x,y)=-xy,\,x,y\in\mathbb{R}$. Let $e$ be a unit vector in
$\mathbb{R}$. Thus $\{e_{1},e_{2},e\}$ is a moving frame in $F$
and the connection 1-form (\ref{newform}) defines a covariant
derivation $\widetilde{\nabla}^{h}$:
\begin{equation}\lb{def1}
\widetilde{\nabla}^{h}e_{1}=-\omega_{2}^{1}\otimes
e_{2}+\omega^{1}\otimes e,\;
\widetilde{\nabla}^{h}e_{2}=\omega_{2}^{1}\otimes
e_{1}+\omega^{2}\otimes e,\;
\widetilde{\nabla}^{h}e=\omega^{1}\otimes e_{1}+\omega^{2}\otimes
e_{2},
\end{equation}
which conserves the metric $\widetilde{g}_{h}$. It is easily seen
that when $M^{2}$ is the hyperbolic plane $\mathbb{H}^{2}$,
imbedded in the standard way as a one sheet of two-sheeted
hyperboloid into the pseudoeuclidean space $ \mathbb{E}^{2,1}$,
then $F$ is simply the trivial bundle
$\mathbb{E}^{2,1}\times\mathbb{H}^{2}$. In this case $E$ is the
normal bundle over hyperboloid
$\mathbb{H}^{2}\subset\mathbb{E}^{2,1}$. This explains the
null-curvature for $\widetilde{\nabla}^{h}$ when
$M^{2}=\mathbb{H}^{2}$, because $\widetilde{\nabla}^{h}$ on
$\mathbb{E}^{2,1}\times\mathbb{H}^{2}$ is the restriction of the
flat Levi-Civita connection on
$T\mathbb{E}^{2,1}=\mathbb{E}^{2,1}\times\mathbb{E}^{2,1}$.

To be sure that the connection $\widetilde{\nabla}^{h}$ is
well-defined on the whole bundle $F$ for the general metric $g$ on
$M^{2}$ we can rewrite $\widetilde{\nabla}^{h}$ as follows. Let
$\xi+f e=\xi^{1}e_{1}+\xi^{2}e_{2}+fe$ is a direct expansion of an
arbitrary section of $F$ over $U$, where $f$ is a smooth function
on $M^{2}$ and $\xi$ is a sections of $TM^{2}$. Then due to
(\ref{def1}) we obtain:
\begin{align}\lb{defSecH}
\widetilde{\nabla}^{h}_{X}(\xi+fe)&=X(\xi^{1})e_{1}+X(\xi^{2})e_{2}
-\xi^{1}\omega_{2}^{1}(X)e_{2}+\xi^{2}\omega_{2}^{1}(X)e_{1}+f(\omega^{1}(X)e_{1}+\omega^{2}(X)e_{2})
\nonumber\\&+ (X(f)+\xi^{1}\omega^{1}(X)+\xi^{2}\omega^{2}(X))e=
\nabla_{X}\xi+fX+\left(X(f)+g(X,\xi)\right)e,
\end{align}
where $X$ is a vector field on $M^{2}$. It is obvious that this
formula gives the definition for $\widetilde{\nabla}^{h}$ on the
whole bundle $F$.

It is possible to change the connection on the bundle $F$ in such
a way that it will be flat iff $g$ is the metric of constant
positive curvature $K=1$. To do this we should write the
connection 1-form $A$ as:
$$
A= \left(\begin{array}{ccc} 0 & \omega_{2}^{1} & \omega^{1} \\
-\omega_{2}^{1} & 0 & \omega^{2} \\ -\omega^{1} & -\omega^{2} & 0
\end{array}\right).
$$
The corresponding derivation will be
\begin{equation}\lb{defSecS}
\widetilde{\nabla}^{s}_{X}(\xi+fe)=\nabla_{X}\xi+fX+\left(X(f)-g(X,\xi)\right)e.
\end{equation}
We see that now $A$ is a $\mathfrak{so}(3)$ valued differential
form and the derivation $\widetilde{\nabla}^{s}$ conserves the
positively defined metric $\widetilde{g}_{s}$ on fibers which is
the direct sum of the metric $g$ on $TM^{2}$ and the metric
$(x,y)=xy,\,x,y\in\mathbb{R}$ on $\mathbb{R}$.

\section{Generalization on higher dimensions}
\label{Generalization}\markright{\ref{Generalization}
Generalization on higher dimensions}

The formulas (\ref{defSecH}) and (\ref{defSecS}) make possible the
immediate generalization of this construction on higher
dimensions. Now $M^{n}$ becomes a $n$-dimensional Riemannian
manifolds and $F$ is the bundle $TM^{n}\oplus E$, where
$E=M^{n}\times\mathbb{R}$. Define the connections
$\widetilde{\nabla}^{s}$ by (\ref{defSecS}) and connection
$\widetilde{\nabla}^{h}$ by the right hand side of
(\ref{defSecH}), where $e$ again is the unit element of the fiber
$\mathbb{R}$ of $F$. The definitions for metrics
$\widetilde{g}_{h}$ and $\widetilde{g}_{s}$ are the same as in the
previous section.

It is well-known that the Riemannian tensor $R$ on a manifold with
constant sectional curvature $K$ satisfies the following relation:
$$
R(X,Y)Z\equiv\nabla_{X}\nabla_{Y}-\nabla_{Y}\nabla_{X}-\nabla_{[X,Y]}=K(g(Y,Z)X-g(X,Z)Y).
$$
We denote such tensor as $R_{K}$. Let us calculate the curvature
tensor $R^{h,s}$ \cite{KoNo1}, corresponding to the connections
$\widetilde{\nabla}^{h}$ and $\widetilde{\nabla}^{s}$ on $F$:
$$
R^{h,s}(X,Y)\widetilde{\xi}\equiv\widetilde{\nabla}^{h,s}_{X}\widetilde{\nabla}^{h,s}_{Y}\widetilde{\xi}-
\widetilde{\nabla}^{h,s}_{Y}\widetilde{\nabla}^{h,s}_{X}\widetilde{\xi}-\widetilde{\nabla}^{h,s}_{[X,Y]}\widetilde{\xi},
$$
where $\widetilde{\xi}=\xi+fe$ is a section of $F$. From
(\ref{defSecS}) we obtain:
\begin{align*}
\widetilde{\nabla}^{s}_{X}\widetilde{\nabla}^{s}_{Y}\left(\xi+fe\right)&=\nabla_{X}\nabla_{Y}\xi+\nabla_{X}(fY)+
\left(Y(f)-g(Y,\xi)\right)X
\\&+\left(X\left(Y(f)-g(Y,\xi)\right)-g(X,\nabla_{Y}\xi+fY)\right)e,
\end{align*}
so
\begin{align*}
R^{s}(X,Y)\widetilde{\xi}&=\nabla_{X}\nabla_{Y}\xi-\nabla_{Y}\nabla_{X}\xi+\nabla_{X}(fY)-\nabla_{Y}(fX)
+Y(f)X-X(f)Y-g(Y,\xi)X\\&+g(X,\xi)Y+\left\{(X\circ Y(f)-Y\circ
X(f)-X\left(g(Y,\xi)\right)+Y\left(g(X,\xi)\right)-g(X,\nabla_{Y}\xi)\right.\\&+\left.g(Y,\nabla_{X}\xi)\right\}e
-\nabla_{[X,Y]}\xi-f[X,Y]-\left([X,Y](f)-g([X,Y],\xi)\right)e=R(X,Y)\xi\\&+f\left(\nabla_{X}Y-\nabla_{Y}X-[X,Y]\right)
-R_{1}(X,Y)\xi+\left(g(\nabla_{Y}X,\xi)-g(\nabla_{X}Y,\xi)\right.\\&+\left.g([X,Y],\xi)\right)e=R(X,Y)\xi-R_{1}(X,Y)\xi,
\end{align*}
due to the equality
$K_{\nabla}(X,Y)=\nabla_{X}Y-\nabla_{Y}X-[X,Y]=0$ for the torsion
$K_{\nabla}$ of the Levi-Civita connection and the condition
$\nabla_{X}g=0$. Reasoning in the similar way, we obtain:
$$
R^{h}(X,Y)\widetilde{\xi}=R(X,Y)\xi-R_{-1}(X,Y)\xi.
$$
Thus the connection $\widetilde{\nabla}^{h}$ is flat iff $M^{n}$
is a space of the constant sectional curvature $-1$ and the
connection $\widetilde{\nabla}^{s}$ is flat iff $M^{n}$ is a space
of the constant sectional curvature $1$.

We can verify that the connection $\widetilde{\nabla}^{h}$
conserves the metric $\widetilde{g}_{h}$ and the connection
$\widetilde{\nabla}^{s}$ conserves the metric $\widetilde{g}_{s}$.
Indeed, let $\widetilde{\xi_{i}}=\xi_{i}+f_{i}e,\,i=1,2$ be
sections of $F$. Then
$$
\widetilde{\nabla}^{s}_{X}\left(\widetilde{g}_{s}(\widetilde{\xi}_{1},\widetilde{\xi}_{2})\right)=
X\left(g(\xi_{1},\xi_{2})+f_{1}f_{2}\right)=g(\nabla_{X}\xi_{1},\xi_{2})+g(\xi_{2},\nabla_{X}\xi_{1})+
X(f_{1}f_{2}).
$$
On the other hand
\begin{align*}
\widetilde{g}_{s}&\left(\widetilde{\nabla}^{s}_{X}\widetilde{\xi}_{1},\widetilde{\xi}_{2}\right)+
\widetilde{g}_{s}\left(\widetilde{\xi}_{1},\widetilde{\nabla}^{s}_{X}\widetilde{\xi}_{2}\right)=
\widetilde{g}_{s}\left(\nabla_{X}\xi_{1}+f_{1}X+\left(X(f_{1})-g(X,\xi_{1})\right)e,\xi_{2}+f_{2}e\right)
\\&+\widetilde{g}_{s}\left(\xi_{1}+f_{1}e,\nabla_{X}\xi_{2}+f_{2}X+\left(X(f_{2})-g(X,\xi_{2})\right)e\right)=
g(\nabla_{X}\xi_{1}+f_{1}X,\xi_{2})\\&+\left(X(f_{1})-g(X,\xi_{1})\right)f_{2}+
g(\xi_{1},\nabla_{X}\xi_{2}+f_{2}X)+\left(X(f_{2})-g(X,\xi_{2})\right)f_{1}\\
&=g(\nabla_{X}\xi_{1},\xi_{2})+g(\xi_{2},\nabla_{X}\xi_{1})+
X(f_{1})f_{2}+X(f_{2})f_{1}.
\end{align*}
The last two equalities give
$$
\left(\widetilde{\nabla}^{s}_{X}\widetilde{g}_{s}\right)(\widetilde{\xi}_{1},\widetilde{\xi}_{2})\equiv
\widetilde{\nabla}^{s}_{X}\left(\widetilde{g}_{s}(\widetilde{\xi}_{1},\widetilde{\xi}_{2})\right)-
\widetilde{g}_{s}\left(\widetilde{\nabla}^{s}_{X}\widetilde{\xi}_{1},\widetilde{\xi}_{2}\right)-
\widetilde{g}_{s}\left(\widetilde{\xi}_{1},\widetilde{\nabla}^{s}_{X}\widetilde{\xi}_{2}\right)=0.
$$
A similar reasoning gives
$\widetilde{\nabla}^{h}_{X}\widetilde{g}_{h}=0$. The conservation
of this metrics means that the corresponding connection 1-form $A$
is $\mathfrak{so}(n+1)$ valued for $\widetilde{\nabla}^{s}$ and
$\mathfrak{so}(n,1)$ valued for $\widetilde{\nabla}^{h}$ with
respect to orthonormal moving frames.

Let now $M^{n}$ be a simply connected space with constant
sectional curvature $\pm 1$. Considering the standard models for
this space as a submanifold of Euclidean (for $K=1$) or
pseudoeuclidean (for $K=-1$) spaces \cite{KoNo1}, we see that the
bundle $F$ is isomorphic to the trivial bundle
$\mathbb{E}^{n,1}\times \mathbb{H}^{n}$ for $K=-1$ and to
$\mathbb{E}^{n+1}\times \mathbb{S}^{n}$ for $K=1$, where
$\mathbb{E}^{n+1}$ is the $(n+1)$-dimensional Euclidean space and
$\mathbb{E}^{n,1}$ is the $(n+1)$-dimensional pseudoeuclidean
space of signature $(n,1)$. In these cases the connection
$\widetilde{\nabla}^{h,s}$ is the restriction of the flat
Levi-Civita connection for
$T\mathbb{E}^{n,1}=\mathbb{E}^{n,1}\times\mathbb{E}^{n,1}$ onto
$\mathbb{E}^{n,1}\times \mathbb{H}^{n}$ or of the flat Levi-Civita
connection for
$T\mathbb{E}^{n+1}=\mathbb{E}^{n+1}\times\mathbb{E}^{n+1}$ onto
$\mathbb{E}^{n+1}\times \mathbb{S}^{n}$.

\section{Discussion}
\label{Discussion}\markright{\ref{Discussion} Discussion} The
common point of view \cite{Sym}, \cite{AntSym} is that the R.
Sasaki connection is based only on internal geometry of surfaces.
However the connection 1-form (\ref{NewConnect}) possesses the
additional (with respect to internal geometry) matrix structure.
Here we interpreted this additional structure as the trivial
one-dimensional summand in the corresponding vector bundle. In the
case of the constant sectional curvature this summand becomes a
normal bundle of the hypersurfaces. Thus our interpretation means
a "virtual" imbedding of the initial space $M^{n}$ into the space
$\mathbb{E}^{n+1}$ or $\mathbb{E}^{n,1}$. This "virtual" imbedding
becomes actual when $M^{n}$ is a space with constant sectional
curvature $\pm 1$.

In the papers \cite{Terng}, \cite{BeTe} a multi-dimensional
generalization of Sine-Gordon equation $u_{xy}=\sin u$ was given
as an imbedding condition of $M^{n}$ into $\mathbb{E}^{2n-1}$. On
the other hand it is well-known \cite{Sasa} that the condition
$dA+A\wedge A$ for matrix 1-form $A$ given by (\ref{NewConnect})
is equivalent to the Sine-Gordon equation whenever differential
forms $\omega^{1},\omega^{2}$ are parameterized by the function
$u$ in a definite way. The generalization of R. Sasaki connection
given in this paper seems to be quite natural, so it can lead to
another multi-dimensional generalization of the Sine-Gordon
equation.

After finishing the present paper the author has received the
letter of Jack Lee from University of Washington (to whom the
author expresses his deep gratitude), who has pointed to the paper
\cite{Min-Oo} of M.~Min-Oo. In that paper the connections
$\widetilde{\nabla}^{h}$ on $TM\oplus E$ was constructed under the
name {\it hyperbolic Cartan connection} in order to prove the
hyperbolic version of the positive mass theorem. There are no any
links with the theory of integrable partial differential equations
and particularly with the Sasaki's construction. Thus the present
paper establishes the connection between pure geometrical
construction of M.~Min-Oo and the well known construction from the
theory of integrable partial differential equations.

\end{document}